\documentclass[12pt]{article}
\usepackage{geometry}                
\usepackage{graphicx}
\usepackage{amssymb}					
\usepackage{amsmath}					
\usepackage{mathrsfs}
\usepackage{tikz}
\usepackage{epstopdf}
\usepackage{amsthm}
\newtheorem*{theorem*}{Theorem}
\newtheorem*{lemma*}{Lemma}
\newtheorem*{example*}{Example}
\newtheorem*{question*}{Question}
\newtheorem*{conjecture*}{Conjecture}
\newtheorem{theorem}{Theorem}

\newtheorem{corollary}[theorem]{Corollary}

\def\b0{{\bf 0}}
\def\b1{{\bf 1}}

\def\n{\noindent}

\begin{document}
\title{Star arboricity relaxed book thickness of $K_n$ \thanks{December 1, 2024}
}
\author{
Paul C. Kainen\\
 \texttt{kainen@georgetown.edu}}
\date{}
\newcommand{\Addresses}{{
  \bigskip
  \footnotesize

\n

\n

\par\nopagebreak
}}

\maketitle

\abstract{

\n {\it
A book embedding of the complete graph $K_n$ needs $\lceil \frac{n}{2} \rceil$ pages and the page-subgraphs can be chosen to be spanning paths (for $n$ even) and one spanning star for $n$ odd. We show that all page-subgraphs can be chosen to be {\rm star forests} by including one extra {\rm cross-cap} page or two new ordinary pages.
}} \\

\n
\n MSC(2020): 05C10; 05C12; 05C78.

\section{Introduction}

Edge decompositions of a graph are often performed so that the parts have required topological properties.  For instance, the {\it thickness} of a graph is the least number of planar subgraphs in a partition of the edge-set.  The {\it book thickness} of a graph requires that the parts of an edge partition (called pages) are outerplane with respect to a fixed ordering of the vertices along a circle.

One might wish to improve the topological and geometric properties of the page subgraphs (e.g., for the purpose of a simple parallel computational architecture).   To get pages that are vertex-disjoint unions of {\it edges} would require at least as many pages as the maximum degree $\Delta$ of the graph so a planar graph could need arbitrarily many such {\it matching} pages, while its book thickness is at most 4 (\cite{yan, yan2}).   But pages that are vertex-disjoint unions of {\it stars} can be achieved.
Indeed, it is well-known that each outerplane graph is the union of two forests and clearly a forest is the union of two star forests.  

Hakimi, Mitchem and Schmeichel \cite{hakimi} showed that each outerplane graph is actually the union of {\it three} star forests.  Hence, for planar graphs, at most 12 star forest pages suffice.

For the complete graph, we can answer the corresponding question.  The book thickness of $K_n$ is $\lceil \frac{n}{2} \rceil$ (\cite{lto, bk79}) and the result of Hakimi et al \cite{hakimi} shows that a multiplicative factor of 3 suffices.
But this is trivially improved by making each page a single star and $n{-}1$ such pages suffice. We give constructions to show that adding {\it one} new ``cross-cap'' and {\it two} new ordinary pages to the minimum set of pages needed for a book embedding of $K_n$ allows the desired star forest attribute to be achieved by the page subgraphs.

Forgetting the topological condition, one may consider the combinatorial problem of finding the least number $r$ of parts for an edge decomposition of $K_n$ into star forests. But $r = 1 + \lceil \frac{n}{2} \rceil$ (\cite{aoki}, \cite{lin-shyu}) so our bound is best-possible.  

The next section has definitions and results; we conclude with a discussion.

\section{Definitions and main result}

The {\bf book thickness} (or {\bf page number}) of a graph $G=(V,E)$ is the least number of parts in an edge-partition of $G$ such that no two same-color edges cross \cite{bk79}. 
The page number of $K_n$ is $\lceil \frac{n}{2} \rceil$, Ollmann \cite{lto}. Indeed, when $n$ is even, this is achieved by $\frac{n}{2}$ Hamiltonian paths; when $n$ is odd, add a star centered at the new vertex - see Harary \cite[pp 91--92, Fig. 9.9]{harary}.  
The {\bf arboricity} $\Gamma$ of $G$ is the least number of forests into which $E$ can be decomposed; $\Gamma(K_n) = \lceil \frac{n}{2} \rceil$.

A {\bf $k$-book} is the regular CW-complex obtained by attaching a family of $k \geq 1$ closed topological disks $D$ to the circle $S^1$, \cite{bk79}.  A {\bf generalized book} (Overbay \cite{so-thesis})  allows the attachment of the complement of an open topological disk in any surface $S$ (orientable or not).  If $S$ is not the sphere, the complement of an open disk is not homeomorphic to a disk.
These disk-complements are called {\bf pages}.  A {\bf cross-cap page} is the complement of an open disk in the projective plane. 

A topological embedding of $G=(V,E)$ in a book or generalized book is called a {\bf (generalized) book embedding} if the vertices are in the bounding circle and each edge is contained in one and only one of the pages.  For a {\bf strict} book embedding (pages are disks), this is equivalent to an outerplane drawing of $G$ and an edge-partition such that each part is crossing-free. A generalized book embedding is {\bf relaxed} if it has a single cross-cap page.

We shall call the triple ($G$, $\omega$, $c$) a  book embedding, where $\omega$ is the cyclic order of $V$  and $c$ is the 
edge partition by pages which may not be a {\it proper edge-coloring} \cite[p 133]{harary} as edges in the same part can share a common endpoint. For a generalized book embedding, one must also identify the page topologies.


The {\bf star arboricity} $sa(K_n)$ of $K_n$ is the least number of induced star forests in any decomposition of $E(K_n)$.  and for $n \geq 4$, we now show that one has $sa(K_n) \geq 1 + \lceil \frac{n}{2} \rceil$.
For suppose that $K_n = F_1 \cup \cdots \cup F_s$ is any edge partition into $sa(K_n)$ spanning forests.  Let $n = 2r$. Then $|E(F_j)| = 2r - \pi_j$, where $\pi_j$ is the number of connected components in $F_j$.  Therefore, one has
\begin{equation}
|E(K_n)| = 2r^2 - r = \sum_{j=1}^s |E(F_j)| \leq s(2r - 1).
\end{equation}
Hence, $s \geq r$, and $s = r$ if and only if each forest has a single component.  If the forests are star forests, then they consist of a single star, but for $n$ odd or even, one needs $n-1$ stars to cover the edges of $K_n$.
 This proves the lower bound in the following result for $n$ even. We continue with the lower bound on star arboricity for $n$ odd.
\begin{theorem}[Aoki \cite{aoki}, Lin \& Shyu \cite{lin-shyu}]
If $n \geq 4$, then $sa(K_n) = 1 + \lceil \frac{n}{2} \rceil$.
\label{th:sa}
\end{theorem}
\begin{proof}
Let $n = 2r+1$.  Then
$|E(K_n)| = 2r^2 + r \leq s(2r)$
so $s \geq r+1 = \lceil \frac{n}{2} \rceil$. But equality again
forces the forests to have one component.
Thus, we have $sa(K_n) \geq r+2 = 1 + \lceil \frac{n}{2} \rceil$.  The upper bound is implictly shown below.
\end{proof}

Let {\bf star arboricity relaxed book thickness} $sarbt(K_n)$ be the least number of pages in a relaxed book embedding of $K_n$ with star forest pages.  We now show that the lower bound on star arboricity is achievable subject to a relaxed book embedding proving equallity in both theorems.

\begin{theorem}
If $n \geq 4$, then $sarbt(K_n) = 1+\lceil \frac{n}{2} \rceil$. 
\label{th:relax}
\end{theorem}
\begin{proof}
By adding a spanning star centered at the new vertex, the upper bound for $n$ odd follows from $n$ even.  Suppose $n = 2r \geq 4$ and put the $n$ vertices $1, 2, \ldots, 2r$ around a circle in counterclockwise order as the $n$ roots of unity.  The cross-cap page accommodates the $r$ edges that join the antipodal vertex pairs.  For each antipodal pair $(i, \,r+i)$, $1 \leq i \leq r$, take the stars $\{ij: i+1 \leq j \leq i+r-1\}$ and $\{[r+i,j]: r+i+1 \leq j \leq i-1\}$.  This gives $r$ star forest pages where each page has two stars, each with $r-1$ edges, and $r$ edges in the cross-cap page.  The total number of edges is $2r^2 - r$ as required.
\end{proof}

The {\bf star arboricity book thickness} $sabt(K_n)$ is the least number of star forest pages in a {\it strict} book embedding of $K_n$. 

\begin{theorem}
If $n \geq 4$, then $1+\lceil \frac{n}{2} \rceil \leq sabt(K_n) \leq 2+\lceil \frac{n}{2} \rceil$.
\label{th:sabt}
\end{theorem}
\begin{proof}
Theorem \ref{th:sa} gives the lower bound. As in the previous proof, the odd case follows from the even case.
For the upper bound on $sabt(K_{2r})$, put the edge $[i, i+r]$ into the star at $i$ and replace the star at $i+r$ with a star at $1+i+r$. The two disjoint stars have $r$ and $r-2$ edges, respectively, and give one of $r = \lceil \frac{n}{2} \rceil$ star forest pages, while the star with $r-1$ edges at $r+1$ and the edge $[r+2, 2r]$ provide two additional star forest pages (one very ``sparse'' \cite{kjo}).  This is a star forest strict book embedding with $\lceil \frac{n}{2} \rceil+2$ pages. As a check, $r(2r-2) + r = {2r \choose 2}$.
\end{proof}

By \cite{bk79}, for any graph $G$ with at least 4 vertices, $bt(G) \geq \frac{|E(G)| - |V(G)|}{|V(G)|-3}$.
Thus, for the {\bf octahedron} $O_r := K_{2r} - rK_2$, $bt(O_r) \geq \lceil \frac{2r^2 -2r - 2r}{2r-3} \rceil = r +\lceil \frac{-r}{2r-3} \rceil$. But $r < 2r-3$ for $r \geq 4$, so $bt(O_r) \geq r$ for $r \geq 4$.  The
proof of Theorem \ref{th:relax} shows this number is achievable by a strict book embedding with $r$ {\it star forest} pages. 
\begin{corollary}
For $r \geq 4$, $sabt(O_r) = r$.
\label{co:octahedron}
\end{corollary}
The octahedron is the $r{-}1$-st power of the cycle $C_{2r}$, where the $k$-th power $G^k$ is the supergraph of $G$ which adds an edge to $G$ for each pair of vertices with $G$-distance in $\{2,\ldots,k\}$ \cite[p 14]{harary}. In \cite{kt-BICA}, we showed that if $n {\geq} 3$ and $1 {\leq} k {<} \frac{n}{2}$, then $sabt(C_n^k) \geq k+1$ and equality holds if $k+1$ divides $n$.  Taking $n=2r$ and $k=r-1$, this gives a different proof for Corollary \ref{co:octahedron}.

Finally, we observe that for $n \leq 5$, $sa(K_n) = sarbt(K_n) = sabt(K_n) = n-1$.

\section{Discussion}

Improvements in geometric and topological properties for the subgraphs in edge decompositions of $K_n$ have been shown attainable by a small {\it additive} increment in the number of parts.  Other examples include the following.

Cycle powers \cite{kt-BICA}. When $k(k+1)$ divides $n$, 
$sabt(C_n^k) = 1 + bt(C_n^k)$; i.e., a single extra page permits star forest pages.  Circulant graphs \cite{alam-2021, kjo, OJK, YSL}.  
For circulant graphs with at most two jumps and for the majority of circulants with longest jump of length 3,
it is possible to find a book embedding whose pages are matchings with at most one additional page above the minimum $\Delta$; for bipartite circulants, there is no increment.  Cartesian products of cycles and complete graphs \cite{pck-BICA, SLL} also have efficient matching book embeddings.

Another example of improvement is the phenomenon of {\it sparse pages} for matching book embeddings: deletion of a small set of edges (less than 6) reduces page number by 1 for  graphs of arbitrarily large size \cite{kjo}.  This applies also to star forest pages;
proof of Theorem \ref{th:sabt} shows that $stabt(K_n - e) = 1 + \lceil \frac{n}{2} \rceil$.

It seems natural to regard page-count as cost but the appropriate ``price'' is not clear for a cross-cap or other non-planar page, nor is it obvious how to value the achieved properties for the pages.
There are also scientific and engineering questions regarding the fabrication of non-planar pages.
Perhaps a cross-cap page might be implemented using a {\it photonic} layer as free-space light beams do not interact with one another; see, e.g., \cite{luo, reynolds}.  
Beginning with a generalized book embedding of $K_n$ on a single sufficiently complex page, one could increase the number of pages to get desired properties for page-graphs.  




Regarding the question of the star forest book thickness of planar graphs, three pages suffice if $G$ is also cubic bipartite \cite{alam-2021, ko} as a matching is a star forest.  How many star forest pages are needed for the icosahedron?  

All of the above can be investigated for other types of geometric constraints.  For instance, ``book'' can be replaced with ``geometric'' meaning that graph drawings are in the plane with vertices no longer confined to being on the unit circle (or other convex boundary) but keeping the requirement that edges are straight-line segments. The largest possible {\it geometric star arboricity} of a planar graph is 5 as \cite{hakimi} gives the upper bound while \cite{aa} shows it is achieved.  But geometric star arboricity for most graphs is unknown currently.

\end{document}